\documentclass{elsart}

\usepackage{amssymb}
\usepackage{amsmath}
\numberwithin{equation}{section}

\begin{document}

\begin{frontmatter}

% Title, authors and addresses

% use the thanksref command within \title, \author or \address for footnotes;
% use the corauthref command within \author for corresponding author footnotes;
% use the ead command for the email address,
% and the form \ead[url] for the home page:
% \title{Title\thanksref{label1}}
% \thanks[label1]{}
% \author{Name\corauthref{cor1}\thanksref{label2}}
% \ead{email address}
% \ead[url]{home page}
% \thanks[label2]{}
% \corauth[cor1]{}
% \address{Address\thanksref{label3}}
% \thanks[label3]{}

\title{On the asymptotically simplicity of periodic eigenvalues and Titchmarsh's formula }

% use optional labels to link authors explicitly to addresses:
% \author[label1,label2]{}
% \address[label1]{}
% \address[label2]{}

\author{Alp Arslan K\i ra\c{c}}
\ead{aakirac@pau.edu.tr}
\address{Department of Mathematics, Faculty of Arts and Sciences, Pamukkale
University, 20070, Denizli, Turkey}

\begin{abstract}
We consider Sturm-Liouville equation $y^{\prime\prime}+(\lambda-q)y=0$ where $q\in L^{1}[0,a]$.
We obtain various conditions on the Fourier coefficients of $q$ such that the periodic eigenvalues having the form given by Titchmarsh are asymptotically simple. Under these conditions, we give some asymptotic estimates for the spectral gaps.
\end{abstract}

\end{frontmatter}

\section{Introduction}
In this study we consider Hill's equation
\begin{equation}  \label{1}
y^{\prime\prime}+(\lambda-q(x))y=0,
\end{equation}
where $q(x)$ is a real-valued summable function on $[0,a]$ for $a>0$ extended to the real line by periodicity. Let
the points $\lambda_{2m+1}, \lambda_{2m+2}$ denote the periodic
eigenvalues of  (\ref{1}) on the interval $[0, a]$ with the
periodic boundary conditions
\begin{equation}\label{per}
y(0)=y(a),\qquad y^{\prime}(0)=y^{\prime}(a).
\end{equation} The number of periodic eigenvalues is countably infinite
 and the eigenvalues form a monotone
increasing sequence with a single accumulation point at $\infty$. See basics and further references in
\cite{Eastham,Coddington:Levinson,Titchmarsh-II,Fulton}.

In the classical investigations, the order of the asymptotic
estimates for the two eigenvalues $\lambda_{2m+1}, \lambda_{2m+2}$
is closely related to the order of smoothness of the potential
$q$. We mention in particular
\cite{Eastham,Titchmarsh-II,Hochstadt:Sturm-Eigenvalues,Naimark}
and the latest results in \cite{Eastham:Eiegen}. In
\cite{Eastham}, Theorem 4.2.4, if, for $r\geq 1$, $q$ has an absolutely continuous
$(r-1)$st derivative on $(-\infty,\infty)$, the asymptotic
estimates for the $\lambda_{2m+1}, \lambda_{2m+2}$ as
$m\rightarrow\infty$ have an error term of order $o(m^{-r})$.

In the case when $r=0$, $q$ is assumed to be piecewise continuous with period
$a$. The presented proof of Theorem 4.2.4 in \cite{Eastham} is
due to Hochstadt \cite{Hochstadt:Sturm-Eigenvalues} by using the
Pr\"{u}fer transformation. Then Titchmarsh
\cite{Titchmarsh-II} (see Section 21.5) refined the $o$-term in
\cite{Eastham}, which leads to the following formulas for the
$\lambda_{2m+1}, \lambda_{2m+2}$, depending on the Floquet theory
(see \cite{Eastham} ,Chapter 2),

\begin{equation}\label{asy1}
{\lambda_{2m+2}\quad\atop \lambda_{2m+1}}\!\!\!\!
=\frac{(2m+2)^{2}\pi^{2}}{a^{2}}\pm |c_{2m+2}|+O(m^{-1/2})
\end{equation}
as $m\rightarrow\infty$, where
\[
c_{2m+2}=:(q\,,e^{i2(2m+2)\pi
x/a})=:\frac{1}{a}\int_{0}^{a}q(x)\,e^{-i2(2m+2)\pi
x/a}\,dx
\]
is the Fourier coefficient of $q$ and without loss of generality we always suppose that
$c_{0}=0$.
Note that the $O$-term in (\ref{asy1}) was essentially given as
$O(m^{-1})$, but Eastham \cite{Eastham} explained, Theorem
4.4.1, that Titchmarsh's proof produces only the $O$-term as in
(\ref{asy1}). For the first time, by using the Pr\"{u}fer transformation,
Brown and Eastham \cite{Eastham:Eiegen} proved, Theorem 2.2, that if $q$
is locally integrable on $(-\infty,\infty)$, then the Titchmarsh's formula (\ref{asy1}) can indeed be
improved to the $O$-term $O(m^{-1}).$ The present work was stimulated by the papers \cite{Eastham:Eiegen,Veliev:Shkalikov,Harris:form}.

By using a perturbation method developed in
\cite{Veliev:Shkalikov,DERNEK:VELIEV,Veliev:Kiraç} we obtain that, under the hypotheses of Theorem \ref{mainthm} (see also (\ref{c1}) and (\ref{c2})),
the $O$-term in (\ref{asy1}) can indeed be improved to
the $O$-term $O(\rho(m)\,m^{-1})$ which also implies $o(m^{-1})$ and all the periodic eigenvalues are asymptotically simple. In Section \ref{remarksec}, using these estimates, the widths of instability intervals are given with the isolated Fourier coefficients of $q$.

Now, let us formulate the subsequent form of the Riemann-Lebesgue lemma. Since the proof of lemma repeats the arguments of the Lemma 6 in \cite{Harris:form}, we omit its proof.
\begin{lem}\label{l1}
 If $q\in L^1[0,a]$ then $\displaystyle\int_{0}^{\,x}q(t)\,e^{i2m\pi t/a}dt\rightarrow 0$ as $|m|\rightarrow\infty$ uniformly in $x$.
\end{lem}

Using this, we define
\begin{equation}\label{m-8}
\rho(m)=:\sup_{0\leq x\leq a} \left| \int_{0}^{x}q(t)\,e^{\mp i2(2m+2)\pi t/a}dt\right|
\end{equation}
and then we get $\rho(m)\rightarrow 0$ as $m\rightarrow \infty$.

In this paper, we prove the following main result:

\begin{thm}\label{mainthm}
Let $q$ be locally integrable on $(-\infty,\infty)$ and assume that the condition
\begin{equation}\label{maincon}
\lim_{m\rightarrow\infty}\frac{\rho(m)}{m\,c_{2m+2}}=0
\end{equation}
holds. Then the large periodic eigenvalues are simple and (\ref{asy1}) holds with the improved $O$-term $O(\rho(m)\,m^{-1})$,
where $\rho(m)$, defined in (\ref{m-8}), is an order of the Fourier coefficient of $q$.
\end{thm}
Clearly, in Theorem \ref{mainthm}, if, instead of (\ref{maincon}), we assume that either the condition \begin{equation}\label{c1}
c_{2m+2}\sim\rho(m),
\end{equation}
where the notation $a_{m}\sim b_{m}$ means that there exist constants $c_{1}$, $c_{2}$ such that $0<c_{1}<c_{2}$ and $c_{1}<|a_{m}/b_{m}|<c_{2}$ for all sufficiently large $m$,
or the condition
\begin{equation}\label{c2}
|c_{2m+2}|>\varepsilon m^{-1}\quad\textrm{with some $\varepsilon>0$}
\end{equation}
holds, then the assertion of Theorem \ref{mainthm} remains valid.

It easily follows from \cite {Eastham}, Theorem 4.2.3, that the
periodic eigenvalues $\lambda_{2m+1}, \lambda_{2m+2}$ are
asymptotically located in pairs, satisfying the following
asymptotic estimate
\begin{equation}\label{asy2}
    \lambda_{2m+1}=\lambda_{2m+2}+o(1)=(2m+2)^{2}\pi^{2}/a^{2}+o(1),
\end{equation}
for sufficiently large $m$. This
estimate implies that the pair of the eigenvalues
$\{\lambda_{2m+1}, \lambda_{2m+2}\}$ is close to the number
$(2m+2)^{2}\pi^{2}/a^{2}$ and isolated from the remaining
spectrum of the problem by a distance of size $m$. In particular, by using (\ref{asy2}),
the following inequality holds
\begin{equation}  \label{dist1}
|\lambda_{m,j}-\frac{(2(m-k)+2)^{2}\pi^{2}}{a^{2}}|>a^{-2}|k||(2m+2)-k|>C\,m,
\end{equation}
for all $k\neq 0,(2m+2)$, $j=1,2$ and $k\in \mathbb{Z}$, where, here and in subsequent relations, we denote $\lambda_{2m+1}$ and $\lambda_{2m+2}$ by $\lambda_{m,1}$ and $\lambda_{m,2}$ respectively for sufficiently large $m$ and $C
$ is positive constant whose
exact value is not essential. For $q=0$,
$\{e^{-i(2m+2)\pi x/a}, e^{i(2m+2)\pi x/a}\}$ is a basis of the
eigenspace corresponding to the eigenvalue
$(2m+2)^{2}\pi^{2}/a^{2}\neq0$ of the problem
(\ref{1})-(\ref{per}) on $[0,a]$.

\section{Preliminaries}
Let us consider the following relation, for sufficiently large
$m$, in order to obtain the values of periodic eigenvalues
$\lambda_{2m+1}, \lambda_{2m+2}$ corresponding to the normalized
eigenfunctions $\Psi_{m,1}(x),\Psi_{m,2}(x)$:
\begin{equation}  \label{m1}
\Lambda_{m,j,m-k}(\Psi_{m,j},e^{i(2(m-k)+2)\pi
x/a})=(q\,\Psi_{m,j},e^{i(2(m-k)+2)\pi x/a}),
\end{equation}
where $\Lambda_{m,j,m-k}=(\lambda_{m,j}-(2(m-k)+2)^{2}\pi^{2}/a^{2})$, $j=1,2.$ 
This relation can be obtained from the equation (\ref{1}), considering $\Psi_{m,j}(x)$ instead of
$y$ and multiplying both sides by $e^{i(2(m-k)+2)\pi x/a}$. Moreover, to iterate (\ref{m1}) we
use the following relations (e.g., see Lemma 1 in \cite{Melda.O}),
\begin{equation}  \label{m2}
(q\,\Psi_{m,j},e^{i (2m+2)\pi
x/a})=\sum_{m_{1}=-\infty}^{\infty}c_{m_1}(\Psi_{m,j},e^{i
(2(m-m_{1})+2)\pi x/a}),
\end{equation}
\begin{equation}  \label{m3}
|(q\,\Psi_{m,j},e^{i (2(m-m_{1})+2)\pi
x/a})|< 3M
\end{equation}
for all sufficiently large $m$, where $m_{1}\in \mathbb{Z}$, $j=1,2$ and $M=\sup_{m\in \mathbb{Z}}|c_{m}|$.

Now, using (\ref{m2}) in (\ref{m1}) for $k=0$
and then isolating the terms with indices $m_{1}=0,(2m+2)$,  we get
\[
 \Lambda_{m,j,m}(\Psi_{m,j},e^{i(2m+2)\pi
x/a})=c_{2m+2}(\Psi_{m,j},e^{-i (2m+2)\pi
x/a})+
\]
\begin{equation}\label{m4}
\qquad\qquad\qquad\qquad\sum_{m_{1}\neq 0,(2m+2)}c_{m_1}(\Psi_{m,j},e^{i
(2(m-m_{1})+2)\pi x/a})
\end{equation}
by the assumption $c_{0}=0$. After using (\ref{m1}) for $k=m_{1}$ in (\ref{m4}), again using (\ref{m2}) with a suitable indices we get the following relation
\begin{equation}\label{m412}
[\Lambda_{m,j,m}- a(\lambda_{m,j})]u_{m,j}=[c_{2m+2}+b(\lambda_{m,j})]v_{m,j}+R(m),
\end{equation}
where $j=1,2,$
\begin{equation}\label{uv}
u_{m,j}=(\Psi_{m,j},e^{i(2m+2)\pi x/a}),\quad v_{m,j}=(\Psi_{m,j},e^{-i(2m+2)\pi x/a}),
\end{equation}
\[
a(\lambda_{m,j})=\sum_{m_{1}}\frac{c_{m_{1}}c_{-m_{1}}}{\Lambda_{m,j,m-m_{1}}},\quad
b(\lambda_{m,j})=\sum_{m_{1}}\frac{c_{m_{1}}c_{2m+2-m_{1}}}{\Lambda_{m,j,m-m_{1}}},
\]
\begin{equation}\label{R}
  R(m)=\sum_{m_{1},m_{2}}\frac{c_{m_{1}}c_{m_{2}}(q\,\Psi_{m,j},e^{i(2(m-m_{1}-m_{2})+2)\pi x/a})}{\Lambda_{m,j,m-m_{1}}\,\Lambda_{m,j,m-m_{1}-m_{2}}}.
\end{equation}
The sums in these formulas are taken over all integers $m_{1},m_{2}$ such that $m_{1},m_{2},m_{1}+m_{2}\neq 0$ and  $m_{1},m_{1}+m_{2}\neq 2m+2.$

Similarly, by considering the other eigenfunction $e^{-i(2m+2)\pi
x/a}$ corresponding to the eigenvalue $(2m+2)^{2}\pi^{2}/a^{2}$ of
the problem (\ref{1})-(\ref{per}) for $q=0$, one can easily obtain
the following relation
\begin{equation}\label{m413}
[\Lambda_{m,j,m}- a'(\lambda_{m,j})]v_{m,j}=[c_{-2m-2}+b'(\lambda_{m,j})]u_{m,j}+R'(m).
\end{equation}
where \[
a'(\lambda_{m,j})=\sum_{m_{1}}\frac{c_{m_{1}}c_{-m_{1}}}{\Lambda_{m,j,m+m_{1}}},\quad
b'(\lambda_{m,j})=\sum_{m_{1}}\frac{c_{m_{1}}c_{-2m-2-m_{1}}}{\Lambda_{m,j,m+m_{1}}},
\]
\begin{equation}\label{R1}
  R'(m)=\sum_{m_{1},m_{2}}\frac{c_{m_{1}}c_{m_{2}}(q\,\Psi_{m,j},e^{i(2(m+m_{1}+m_{2})+2)\pi x/a})}{\Lambda_{m,j,m+m_{1}}\,\Lambda_{m,j,m+m_{1}+m_{2}}},
\end{equation}
and the sums in these formulas are taken over all integers $m_{1},m_{2}$ such that $m_{1},m_{2},m_{1}+m_{2}\neq 0$ and  $m_{1},m_{1}+m_{2}\neq -2m-2.$

By using (\ref{dist1}), (\ref{m1}) and (\ref{m3}), we get (see \cite[Theorem 2]{Melda.O})
\[
\sum_{k\in \mathbb{Z};\,k\neq \pm(m+1)}\Big|(\Psi_{m,j},e^{i2k\pi x/a})\Big|^{2}=O\left(\frac{1}{m^{2}}\right)
\]
and by using the equality
\[
  \frac{1}{m_{1}(2m+2-m_{1})}=\frac{1}{2m+2}\left(\frac{1}{m_{1}}+\frac{1}{2m+2-m_{1}}\right),
\]
we can easily be shown the following relation
\begin{equation}\label{main}
\sum_{m_{1}\neq 0,(2m+2)}\frac{1}{|m_{1}||2m+2-m_{1}|}=O\left(\frac{ln|m|}{m}\right).
\end{equation}
Therefore, we obtain that the normalized eigenfunctions  $\Psi_{m,j}(x)$ have an expansion in terms of the orthonormal basis $\{e^{i2k\pi x/a}:k\in \mathbb{Z}\}$ on $[0,a]$ of the form (see also related result in (78) on p. 77 of \cite{Naimark})
\begin{equation}\label{m7}
\Psi_{m,j}(x)=u_{m,j}\,e^{i(2m+2)\pi x/a}+v_{m,j}\,e^{-i(2m+2)\pi x/a}+h_{m}(x),
\end{equation}
where
\begin{equation}\label{m81}
 \!\!(h_{m},e^{\mp i(2m+2)\pi x/a})=0,\, \|h_{m}\|=O(m^{-1}),\, \sup_{x\in[0,a]}|h_{m}(x)|=O\left(\frac{ln|m|}{m}\right)
\end{equation}
\begin{equation}\label{m8}
  |u_{m,j}|^{2}+|v_{m,j}|^{2}=1+O\left(m^{-2}\right).
\end{equation}

\section{Estimates for the eigenvalues}\label{results}
The following estimates play an essential role in the proof of main result of the paper.
\begin{lem}\label{L1}
The eigenvalues $\lambda_{2m+1}, \lambda_{2m+2}$ of the problem (\ref{1})-(\ref{per})  satisfy, for
$m\geq N$,
\begin{equation}\label{lem}
\lambda_{2m+1}, \lambda_{2m+2}=(2m+2)^{2}\pi^{2}/a^{2}+O(\rho(m)),
\end{equation}
where $\rho(m)$ is defined in (\ref{m-8}).
\end{lem}
\begin{pf}
From the relations (\ref{dist1}) and (\ref{main}), one can easily see that
\[
\sum_{m_{1}\neq 0,\pm(2m+2)}\left|\frac{1}{\Lambda_{m,j,m\mp m_{1}}}-\frac{1}{\Lambda_{m,m\mp m_{1}}^{0}}\right|
\]
\begin{equation}\label{dif}
\leq C \frac{|\Lambda_{m,j,m}|}{m}\sum_{m_{1}\neq 0,\pm(2m+2)}\frac{1}{|m_{1}||2m+2\mp m_{1}|} =O\left(\frac{\Lambda_{m,j,m}}{m}\right),
\end{equation}
where $\Lambda_{m,m\mp m_{1}}^{0}=((2m+2)^{2}\pi^{2}/a^{2}-(2(m\mp m_{1})+2)^{2}\pi^{2}/a^{2})$.
Thus, we get
\begin{equation}\label{d1}
    a(\lambda_{m,j})=\frac{a^{2}}{4\pi^{2}}\sum_{m_{1}\neq 0,(2m+2)}\frac{c_{m_{1}}c_{-m_{1}}}{m_{1}(2m+2-m_{1})}+O\left(\frac{\Lambda_{m,j,m}}{m}\right).
\end{equation}
It also follows from \cite{veliev;arþiv} (see Lemma 2) that, in our notations,
\[
\!\!\!\!\!\!\!\!\!a(\lambda_{m,j})=\frac{a^{2}}{2\pi^{2}}\sum_{m_{1}> 0,m_{1}\neq(2m+2)}\frac{c_{m_{1}}c_{-m_{1}}}{(2m+2+m_{1})(2m+2-m_{1})}+O\left(\frac{\Lambda_{m,j,m}}{m}\right)
\]
\begin{equation}\label{d3}
=a^{2}\int_{0}^{a}(G^{+}(x,m)-G^{+}_{0}(m))^{2}\,e^{i2(4m+4)\pi x/a}\,dx+O\left(\frac{\Lambda_{m,j,m}}{m}\right),
\end{equation}
where
\begin{equation}\label{d4}
G^{\pm}_{m_{1}}(m)=:(G^{\pm}(x,m), e^{i2m_{1}\pi x/a})=\frac{c_{m_{1}\pm(2m+2)}}{i2\pi m_{1}}
\end{equation}
for $m_{1}\neq 0$ and $G^{\pm}_{0}(m)=:(G^{\pm}(x,m), 1)$ for $m_{1}=0$ are the Fourier coefficients of the functions
\begin{equation}\label{d2}
G^{\pm}(x,m)=\frac{1}{a}\int_{0}^{x}q(t)\,e^{\mp i2(2m+2)\pi t/a}dt-\frac{1}{a}\,c_{\pm(2m+2)}x
\end{equation}
with respect to the trigonometric
system $\{e^{i2m_{1}\pi x/a}: m_{1}\in\mathbb{Z}\}$
and
\begin{equation}\label{gsum}
G^{\pm}(x,m)-G^{\pm}_{0}(m)=\sum_{m_{1}\neq(2m+2)}\frac{c_{m_{1}}}{i2\pi(m_{1}\mp(2m+2))}\,e^{i2(m_{1}\mp(2m+2))\pi x/a}.
\end{equation}
Now, taking into account the equalities (\ref{m-8}) and (\ref{d2}), we obtain the estimates
\begin{equation}\label{gg}
G^{\pm}(x,m)-G^{\pm}_{0}(m)=G^{\pm}(x,m)-\frac{1}{a}\int_{0}^{a}G^{\pm}(x,m)\, dx=O(\rho(m)),
\end{equation}
\[
G^{\pm}(a,m)=G^{\pm}(0,m)=0.
\]
 Then, these with integration by parts and $q\in L^1[0,a]$ imply that (see (\ref{d3}))
\[
 a(\lambda_{m,j})=\frac{-a^{2}}{i2\pi(4m+4)}\times
\]
\[
\!\!\int_{0}^{a}2(G^{+}(x,m)-G^{+}_{0}(m))(q(x)e^{-i2(2m+2)\pi x/a}-c_{2m+2})e^{i2(4m+4)\pi x/a}dx+O\left(\frac{\Lambda_{m,j,m}}{m}\right)
\]
\begin{equation}\label{s0}
    =O\left(\frac{\rho(m)}{m}\right)+O\left(\frac{\Lambda_{m,j,m}}{m}\right)
\end{equation}
for sufficiently large $m$.% Here and in subsequent relations, we have written $\rho(m)$ for $\rho(2m+2)$.

 Similarly, by virtue of (\ref{dif}) we get
\[\!\!\!\!\!\!\!\!\!\!\!\!\!\!\!\!\!\!\!\!\!\!\!\!\!\!\!\!\!\!\!\!\!\!\!\!\!\!\!\!b(\lambda_{m,j})=\frac{a^{2}}{4\pi^{2}}\sum_{m_{1}\neq 0,(2m+2)}\frac{c_{m_{1}}c_{2m+2-m_{1}}}{m_{1}(2m+2-m_{1})}+O\left(\frac{\Lambda_{m,j,m}}{m}\right)\]
\[
\!\!\!\!\!\!\!\!\!\!\!\!\!\!\!\!\!\!\!\!\!\!\!\!=-a^{2}\int_{0}^{a}(Q(x)-Q_{0})^{2}\,e^{-i2(2m+2)\pi x/a}dx+O\left(\frac{\Lambda_{m,j,m}}{m}\right)
\]
\begin{equation}\label{I0}
\qquad\;\;=\frac{-a^{2}}{i2\pi(2m+2)}\int_{0}^{a}2(Q(x)-Q_{0})\,q(x)\,e^{-i2(2m+2)\pi x/a}dx+O\left(\frac{\Lambda_{m,j,m}}{m}\right),
\end{equation}
where  $Q(x)=a^{-1}\displaystyle\int_{0}^{x}q(t)\, dt,\quad Q_{m_{1}}=:(Q(x),e^{i2m_{1}\pi x/a})=\frac{c_{m_{1}}}{i2\pi m_{1}} \quad \textrm{if}\,\, m_{1}\neq 0$,
\begin{equation}\label{Q0}
Q(x)-Q_{0}=\sum_{m_{1}\neq 0}Q_{m_{1}}\,e^{i2m_{1}\pi x/a}.
\end{equation}
Therefore, again using integration by parts the integral in (\ref{I0}), $Q(a)=c_{0}=0$  and (\ref{m-8}), we obtain
\begin{equation}\label{b1}
    b(\lambda_{m,j})=O\left(\frac{\rho(m)}{m}\right)+O\left(\frac{\Lambda_{m,j,m}}{m}\right).
\end{equation}
Also
\begin{equation}\label{b1'}
    b'(\lambda_{m,j})=O\left(\frac{\rho(m)}{m}\right)+O\left(\frac{\Lambda_{m,j,m}}{m}\right).
\end{equation}
Let us prove that $R(m)=O\left(\rho(m)\right)$ (see (\ref{R})). First, since $q\in L^{1}[0,a]$ and is different from $0$, there exists $x\in[0,a]$ such that
\begin{equation}\label{qint}
    \int_{0}^{x}q(t)\,dt\neq 0
\end{equation}
and the integral (\ref{qint}) is finite for all $x\in [0,a].$
Therefore, multiplying  the integrand of (\ref{qint}) by  $e^{-i2(2m+2)\pi x/a}e^{i2(2m+2)\pi x/a}$ and using integration by parts, we have
\begin{equation*}
    \sup_{0\leq x\leq a}\left|\int_{0}^{x}q(t)\,dt\right|\leq C\left(\rho(m)+m\rho(m)\right)\leq Cm\rho(m)
\end{equation*}
which implies that, for sufficiently large $m$,
\begin{equation}\label{rhoineq}
    \rho(m)>C\,m^{-1}
\end{equation}
with some $C>0$. Then using (\ref{dist1}), (\ref{m3}) and the relation (\ref{main}),
we obtain
\begin{equation*}
|R(m)|\leq C\frac{(ln|m|)^{2}}{m^{2}}=O\left(\rho(m)\right).
\end{equation*}
Similarly, $R'(m)=O\left(\rho(m)\right).$

It may readily be seen by a change the sign of summation indices in the relation for $a'(\lambda_{m,j})$ (see (\ref{m412}), (\ref{m413})) that $a(\lambda_{m,j})=a'(\lambda_{m,j}).$
On the other hand, the formula (\ref{m8}) gives that either $|u_{m,j}|>1/2$ or $|v_{m,j}|>1/2$ for large $m.$ Thus If $|u_{m,j}|>1/2$, then using (\ref{m412}), (\ref{s0}) and (\ref{b1}) with $R(m)=O\left(\rho(m)\right)$ we get
 \begin{equation*}
\Lambda_{m,j,m}\left(1+O\left(\frac{1}{m}\right)\right)=c_{2m+2}\frac{v_{m,j}}{u_{m,j}}+O\left(\rho(m)\right).
\end{equation*}
This with (\ref{m-8}) implies $\Lambda_{m,j,m}=O\left(\rho(m)\right).$ If $|v_{m,j}|>1/2$ then, again by (\ref{m413}), (\ref{s0}), (\ref{b1'}) and $R'(m)=O\left(\rho(m)\right)$, we have (\ref{lem}). The lemma is proved.
$\blacksquare$
\end{pf}
\begin{lem}\label{L2}
For all sufficiently large $m$, the series in (\ref{R}) and (\ref{R1}) have the following estimates
\begin{equation}\label{Rlmma}
R(m)=O\left(\rho(m)m^{-1}\right),\qquad R'(m)=O\left(\rho(m)m^{-1}\right).
\end{equation}
\end{lem}
\begin{pf}
First, we prove that $R(m)=O\left(\rho(m)m^{-1}\right)$. Arguing as in the proof of (\ref{d1}) and using Lemma \ref{L1} we get
\begin{equation}\label{Rthm}
R(m)=\sum_{m_{1},m_{2}}\frac{a^{4}2^{-4}\pi^{-4}c_{m_{1}}c_{m_{2}}(q\,\Psi_{m,j},e^{i(2(m-m_{1}-m_{2})+2)\pi x/a})}{m_{1}(2m+2-m_{1})(m_{1}+m_{2})(2m+2-(m_{1}+m_{2}))}+O\left(\frac{\rho(m)}{m^{2}}\right).
\end{equation}
Then, since the eigenfunction $\Psi_{m,j}(x)$ has expansion (\ref{m7}), we have
\begin{equation}\label{Rthm1}
  R(m)=\frac{a^{4}}{(2\pi)^{4}}\left(S(m)u_{m,j}+I(m)v_{m,j}+r(m)\right)+O\left(\frac{\rho(m)}{m^{2}}\right),
\end{equation}
where $S(m)$, $I(m)$ and $r(m)$ are obtained from the series in (\ref{Rthm}) by replacing the numerator of the fraction in the series  with $c_{m_{1}}c_{m_{2}}c_{-m_{1}-m_{2}}$, $c_{m_{1}}c_{m_{2}}c_{2m+2-m_{1}-m_{2}}$ and $c_{m_{1}}c_{m_{2}}(q\,h_{m},e^{i(2(m-m_{1}-m_{2})+2)\pi x/a})$ respectively.

Now to estimate $R(m)$ it is enough to show that
\begin{equation}\label{total}
|S(m)|+|I(m)|+|r(m)|=O\left(\rho(m)m^{-1}\right),
\end{equation}
since $|u_{m,j}|\leq 1$ and $|v_{m,j}|\leq 1$ (see (\ref{uv})) for the normalized eigenfunctions $\Psi_{m,j}(x)$.
 By using the summation variable $m_{2}$ to represent the previous $m_{1}+m_{2}$ in the formula $S(m)$ which is obtained from the series in (\ref{Rthm}), we get
\begin{equation}\label{I1}
  S(m)=\sum_{m_{1},m_{2}}\frac{c_{m_{1}}c_{m_{2}-m_{1}}c_{-m_{2}}}{m_{1}(2m+2-m_{1})m_{2}(2m+2-m_{2})}.
\end{equation}
By using the equality
\[
  \frac{1}{k(2m+2-k)}=\frac{1}{2m+2}\left(\frac{1}{k}+\frac{1}{2m+2-k}\right),
\]
we have
\begin{equation}\label{equal1}
   S(m)=\frac{1}{(2m+2)^{2}}\sum_{j=1}^{4}S_{j},
\end{equation}
where
\[
S_{1}=\sum_{m_{1},m_{2}}\frac{c_{m_{1}}c_{m_{2}-m_{1}}c_{-m_{2}}}{m_{1}m_{2}},\quad S_{2}=\sum_{m_{1},m_{2}}\frac{c_{m_{1}}c_{m_{2}-m_{1}}c_{-m_{2}}}{m_{2}(2m+2-m_{1})},
\]
\[
S_{3}=\sum_{m_{1},m_{2}}\frac{c_{m_{1}}c_{m_{2}-m_{1}}c_{-m_{2}}}{m_{1}(2m+2-m_{2})},\; S_{4}=\sum_{m_{1},m_{2}}\frac{c_{m_{1}}c_{m_{2}-m_{1}}c_{-m_{2}}}{(2m+2-m_{1})(2m+2-m_{2})}.
\]
Now arguing as in the proof of (\ref{b1}) (see definition of (\ref{I0})), we get
\begin{equation}\label{s11}
 S_{1}=-4\pi^{2}\int_{0}^{a}(Q(x)-Q_{0})^{2}q(x)\,dx=O\left(m\rho(m)\right).
\end{equation}
Here, the $O$-term in (\ref{s11}) is obtained from (\ref{rhoineq}), since the integral in (\ref{s11}) is $O(1)$.
Moreover, by (\ref{d4})-(\ref{gg}) and (\ref{Q0}), we have the following relations
\begin{equation}\label{a}
\left.
\begin{array}{ll}
 \displaystyle S_{2}=-4\pi^{2}\int_{0}^{a}(Q(x)-Q_{0})(G^{+}(x,m)-G^{+}_{0}(m))\,q(x)\,e^{i2(2m+2)\pi x/a}dx=O\left(\rho(m)\right), & \\\\
 \displaystyle S_{3}=-4\pi^{2}\int_{0}^{a}(Q(x)-Q_{0})(G^{-}(x,m)-G^{-}_{0}(m))\,q(x)\,e^{-i2(2m+2)\pi x/a}dx=O\left(\rho(m)\right), &  \\\\
 \displaystyle  S_{4}=4\pi^{2}\int_{0}^{a}(G^{+}(x,m)-G^{+}_{0}(m))(G^{-}(x,m)-G^{-}_{0}(m))\,q(x)\,dx=O\left(\rho(m)\right). &
\end{array}
    \right\}
\end{equation}
Thus, in view of (\ref{equal1})-(\ref{a}), we obtain $S(m)=O\left(\rho(m)m^{-1}\right).$
From \cite{Veliev:Shkalikov}, to estimate $I(m)$, using substitutions $k_{1}=m_{1}$, $k_{2}=2m+2-m_{1}-m_{2}$ in the formula $I(m)$
and arguing as in (\ref{I1})-(\ref{equal1}), one can easily obtain that
\begin{equation}\label{equal2}
   I(m)=\frac{1}{(2m+2)^{2}}(I_{1}+2I_{2}+I_{3}),
\end{equation}
where\[
I_{1}=\sum_{m_{1},m_{2}}\frac{c_{m_{1}}c_{m_{2}}c_{2m+2-m_{1}-m_{2}}}{m_{1}m_{2}},\quad I_{2}=\sum_{m_{1},m_{2}}\frac{c_{m_{1}}c_{m_{2}}c_{2m+2-m_{1}-m_{2}}}{m_{2}(2m+2-m_{1})},
\]
\[
I_{3}=\sum_{m_{1},m_{2}}\frac{c_{m_{1}}c_{m_{2}}c_{2m+2-m_{1}-m_{2}}}{(2m+2-m_{1})(2m+2-m_{2})}.
\]
Again by (\ref{d4})-(\ref{gg}), (\ref{Q0}) and the integration by parts only in $I_{1}$, we get the following estimates
\begin{equation}\label{I123}
\left.
\begin{array}{ll}
\displaystyle I_{1}=-4\pi^{2}\int_{0}^{a}(Q(x)-Q_{0})^{2}\,q(x)\,e^{-i2(2m+2)\pi x/a}dx=O\left(\rho(m)\right), &  \\\\
\displaystyle I_{2}=4\pi^{2}\int_{0}^{a}(Q(x)-Q_{0})(G^{+}(x,m)-G^{+}_{0}(m))\,q(x)\,dx=O\left(\rho(m)\right), & \\\\
\displaystyle  I_{3}=-4\pi^{2}\int_{0}^{a}(G^{+}(x,m)-G^{+}_{0}(m))^{2}\,q(x)\,e^{i2(2m+2)\pi x/a}dx=O\left(\rho(m)\right). &
\end{array}
    \right\}
\end{equation}
Thus, from (\ref{equal2}) and (\ref{I123}), $I(m)=O\left(\rho(m)m^{-1}\right).$
Let us prove that
\begin{equation*}
r(m)=\sum_{m_{1},m_{2}}\frac{c_{m_{1}}c_{m_{2}}(q\,h_{m},e^{i(2(m-m_{1}-m_{2})+2)\pi x/a})}{m_{1}(2m+2-m_{1})(m_{1}+m_{2})(2m+2-(m_{1}+m_{2}))}=O\left(\frac{\rho(m)}{m}\right).
\end{equation*}
Now, by using the relation (\ref{main}), together with (\ref{m81}), and taking into account the estimate (\ref{rhoineq}), we have
\[
|r(m)|\leq C\frac{(ln|m|)^{2}}{m^{2}}\sup_{x\in[0,a]}|h_{m}(x)|= \rho(m)\frac{(ln|m|)^{2}}{m}\,O\left(\frac{ln|m|}{m}\right)=O\left(\frac{\rho(m)}{m}\right).
\]
where the positive constant $C$ is independent of $m$. Thus (\ref{total}) holds. Therefore, the first estimate of (\ref{Rlmma}) follows from (\ref{Rthm1}) and (\ref{total}). In the same way,
one can easily obtain the second estimate of (\ref{Rlmma}).
$\blacksquare$
\end{pf}
Now using these lemmas, let us prove the main result of the paper.
\renewenvironment{pf}[1]{{\bfseries Proof of Theorem \ref{mainthm}.#1}}{}

\begin{pf}\, In view of Lemma \ref{L1}, substituting the estimates of
\[a(\lambda_{m,j}),\,a'(\lambda_{m,j}),\,  b(\lambda_{m,j}),\,  b'(\lambda_{m,j}),\, R(m),\, R'(m)=O\left(\rho(m)m^{-1}\right)\]
given by (\ref{s0}), (\ref{b1}), (\ref{b1'}), (\ref{Rlmma})  in the relations (\ref{m412}) and (\ref{m413}), we get the following form of the relations
\begin{equation}\label{son}
\Lambda_{m,j,m}\,u_{m,j}=c_{2m+2}\,v_{m,j}+O(\rho(m)m^{-1}),
\end{equation}
\begin{equation}\label{son'}
\Lambda_{m,j,m}\,v_{m,j}=c_{-2m-2}\,u_{m,j}+O(\rho(m)m^{-1})
\end{equation}
for $j=1,2$. Again by (\ref{m8}), we have, for large $m$, either $|u_{m,j}|>1/2$ or $|v_{m,j}|>1/2$. Without loss of generality we assume that $|u_{m,j}|>1/2$. Then it follows from both (\ref{maincon}), (\ref{son}) and (\ref{maincon}), (\ref{son'}) that
\begin{equation}\label{sameo}
\Lambda_{m,j,m}\sim c_{2m+2},
\end{equation}
where the notation $a_{m}\sim b_{m}$ is defined in (\ref{c1}).
This with (\ref{maincon}), (\ref{son'}) and the assumption $|u_{m,j}|>1/2$ implies that
\begin{equation}\label{vsimilar}
u_{m,j}\sim v_{m,j}\sim 1.
\end{equation}
Thus, first multiplying (\ref{son'}) by $c_{2m+2}$ and then using (\ref{son}) in (\ref{son'}), we get
\[
\Lambda_{m,j,m}(\Lambda_{m,j,m}\,u_{m,j}+O(\rho(m)m^{-1}))=|c_{2m+2}|^{2}\,u_{m,j}+c_{2m+2}O(\rho(m)m^{-1}),
\]
which implies, in view of (\ref{sameo})-(\ref{vsimilar}), the following equations
\begin{equation}\label{final}
\Lambda_{m,j,m}=\pm|c_{2m+2}|+O(\rho(m)m^{-1})
\end{equation}
for $j=1,2$. Arguing as in Lemma 4 of \cite{Veliev:Shkalikov}, let us prove the large periodic eigenvalues are simple.
For large $m$, suppose that there exist two mutually orthogonal eigenfunctions $\Psi_{m,1}(x)$ and $\Psi_{m,2}(x)$ corresponding to $\lambda_{m,1}=\lambda_{m,2}$. Hence, taking into
account the expansion (\ref{m7}) with $\|h_{m}\|=O(m^{-1})$ (see (\ref{m81})) for both the eigenfunctions $\Psi_{m,j}(x)$ and then using their orthogonality, we can choose the eigenfunction $\Psi_{m,j}(x)$ such that either $u_{m,j}=0$ or $v_{m,j}=0$,
which contradicts (\ref{vsimilar}).

Finally, for large $m$, let us prove that each of the simple eigenvalues  in (\ref{final}) corresponds
to only either the lower sign $-$ or the upper sign $+$, not both. In the first case, we assume that both  eigenvalues correspond
to the lower sign $-$. Then by (\ref{son}) and (\ref{final}), we get
\begin{equation}\label{fark1}
(\Lambda_{m,2,m}-\Lambda_{m,1,m})\,u_{m,2}=c_{2m+2}\,v_{m,2}+|c_{2m+2}|\,u_{m,2}+O(\rho(m)m^{-1}),
\end{equation}
\begin{equation}\label{fark2}
(\Lambda_{m,2,m}-\Lambda_{m,1,m})\,u_{m,1}=-c_{2m+2}\,v_{m,1}-|c_{2m+2}|\,u_{m,1}+O(\rho(m)m^{-1}).
\end{equation}

Therefore, using the assumption $(\Lambda_{m,2,m}-\Lambda_{m,1,m})=O(\rho(m)m^{-1})$, multiplying both sides of (\ref{fark1}) and (\ref{fark2}) by $v_{m,1}$ and $v_{m,2}$,
respectively, and then adding them together, we have, in view of (\ref{maincon}),
\begin{equation}\label{uvfark1}
u_{m,2}v_{m,1}-u_{m,1}v_{m,2}=o(1).
\end{equation}
Since the eigenfunctions $\Psi_{m,1}$ and $\overline{\Psi_{m,2}}$  of the self-adjoint problem which belong to
different eigenvalues $\lambda_{m,1}\neq\lambda_{m,2}$ are orthogonal we have, in view of (\ref{m7}),
\[
u_{m,2}v_{m,1}+u_{m,1}v_{m,2}=O(m^{-1}).
\]

This with (\ref{uvfark1}) gives $u_{m,2}v_{m,1}=o(1)$ which contradicts (\ref{vsimilar}). The other case, where both  eigenvalues correspond
to the upper sign $+$, is also impossible. Thus, in
(\ref{final}), we may choose the lower sign $-$ for $j=1$ and the upper sign $+$ for $j=2$. The theorem is proved.
$\blacksquare$
\end{pf}
\renewenvironment{pf}[1]{{\bfseries Proof. #1}}{}

\section{Some remarks}\label{remarksec}

1. Applying the same scheme in the paper to the antiperiodic boundary conditions
\begin{equation*}\label{antiper}
\qquad\qquad\qquad y(0)=-y(a),\qquad y^{\prime}(0)=-y^{\prime}(a),
\end{equation*}
one can easily obtain similar results in Section \ref{results} for the antiperiodic eigenvalues $\mu_{2m}$, $\mu_{2m+1}$ from (\ref{maincon}),(\ref{asy2}), (\ref{dist1}) and (\ref{m1}) by replacing $2m+2$ with $2m+1$. Then we get
$(2m+1)^{2}\pi^{2}/a^{2}$ and $c_{2m+1}$ in (\ref{asy1}) with the $O$-term as $O(\rho(m)\,m^{-1})$.

2. It is well known that the intervals $(\mu_{2m}$, $\mu_{2m+1})$ and $(\lambda_{2m+1}, \lambda_{2m+2})$ are called the instability
intervals of the operator (\ref{1}) generated in $L^{2}(-\infty, \infty)$ by periodicity of $q$. Using Theorem \ref{mainthm}, (\ref{c1}) and (\ref{c2}), the following asymptotic
behaviors of the widths $\ell_{2m+2}$ of these instability intervals are determined for sufficiently large $m$. If either condition (\ref{maincon}) or (\ref{c1}) holds, then we have, respectively,
$$\ell_{2m+2}=c_{2m+2}(1+o(1)),$$
$$\ell_{2m+2}=c_{2m+2}(1+O(m^{-1})),$$
and if the condition (\ref{c2}) holds for sufficiently large $m$, then we have
$$\ell_{2m+2}=c_{2m+2}(1+O(\rho(m)),$$
where $\rho(m)$ is defined in (\ref{m-8}).

\section*{Acknowledgments}
The author is grateful to the anonymous reviewers for their helpful
comments and suggestions.
%\newproof{pot1}{Proof of Theorem \ref{main2}}
%\begin{pot}
%\end{pot}
%\begin{pot1}
%\end{pot1}
 %% The Appendices part is started with the command \appendix;
%% appendix sections are then done as normal sections
%% \appendix

%% \section{}
%% \label{}
%\begin{acknowledgements}
%The author is grateful to the anonymous reviewers for their helpful
%comments and suggestions.
%\end{acknowledgements}
%\bibliographystyle{elsarticle-num}
%\bibliography{article}

\end{document}